\let\ORIlabel\label
\let\ORIrefstepcounter\refstepcounter
\AddToHook{package/hyperref/before}{%
  \let\label\ORIlabel
  \let\refstepcounter\ORIrefstepcounter
}
\documentclass[
onefignum,onetabnum]{siamonline171218}

\usepackage{amsfonts}
\usepackage{amssymb}
\usepackage{amsmath}
\usepackage{mathtools}
\usepackage{enumitem}

\newcommand{\X}{\mathcal X}
\newcommand{\Y}{\mathcal Y}
\newcommand{\Z}{\mathcal Z}
\newcommand{\V}{\mathcal V}

\newcommand{\LL}{\mathcal L}
\newcommand{\dom}{\operatorname{dom}}
\newcommand{\ran}{\operatorname{ran}}
\newcommand{\Real}{\operatorname{Re}}

\newcommand{\Oop}{\mathcal O}
\newcommand{\formh}{\mathfrak h}

\DeclarePairedDelimiter{\norm}{\lVert}{\rVert}
\DeclarePairedDelimiterX{\inner}[2]{\langle}{\rangle}{#1,#2}
\newsiamremark{remark}{Remark}

\headers{Dissipative Hamiltonian representations}{V. Mehrmann and T. Reis}
\title{Dissipative Hamiltonian representations of linear evolution equations}
\author{Volker Mehrmann\thanks{Institut f\"ur Mathematik,
IMoS/MA4-5, Technische Universit\"at Berlin,
Fasanenstr. 89, 10623 Berlin, Germany
(\email{mehrmann@math.tu-berlin.de}).}
\and
Timo Reis\thanks{Fachgebiet Systemtheorie und partielle Differentialgleichungen,
Institut f\"ur Mathematik, Technische Universit\"at Ilmenau,
Postfach 10 05 65, 98684 Ilmenau, Germany
(\email{timo.reis@tu-ilmenau.de}).}}

\ifpdf
\hypersetup{pdftitle={Dissipative Hamiltonian representations of linear evolution equations},
 pdfauthor={V. Mehrmann and T. Reis}}
\fi

\begin{document}
\maketitle

\begin{abstract}
We investigate explicit and implicit dissipative Hamiltonian
representations of generators of strongly continuous semigroups
on real or complex Hilbert spaces. For a prescribed positive definite
selfadjoint energy operator $H$ satisfying
$\dom(A)\subset\dom(H)$, we show that the Lyapunov
dissipation inequality is equivalent to an exact factorization
$A=DH$, where $D$ can be chosen maximally dissipative.
For bounded $H$, the criterion also yields an implicit
representation. We further study energy constructions based
on closed intertwiners and output trajectories, distinguishing
contractive energy-space dynamics from operator factorizations
on the original state space. Examples illustrate the roles
of noncoercivity, unboundedness, and domain compatibility.
\end{abstract}



\section{Introduction}\label{sec:intro}

Dissipative Hamiltonian representations of dynamical systems express the evolution of a system
in terms of an energy operator and a dissipation operator. Besides making
energy balances transparent, these representations provide a structural
starting point for investigating well-posedness, stability, and the properties of related
operator pencils. For finite-dimensional linear differential equations
and differential-algebraic equations, the relationships between various
Hamiltonian, dissipative, and relation-based descriptions have been
studied extensively; see, for instance,
\cite{BeMeXuZw18,GerHaRe21,MMW18,MehrmannVanderSchaft2023,MehS26_ppt}.
Infinite-dimensional counterparts arise in the analysis of partial
differential equations (PDEs) and abstract port-Hamiltonian systems. The formulation of
PDEs as abstract dissipative Hamiltonian differential-algebraic equations
on extended state spaces and the spectral theory of the associated operator pencils have been
investigated in \cite{MehlMehrmannWojtylak2026,ZwaM24};
see also \cite{PhilippReisSchaller2025,Reis2026Dissipation}.
Explicit and implicit port-Hamiltonian formulations of particular
distributed-parameter models have been compared in
\cite{BendimeradHohlMatignonHaineLefevre2024}.
The existing literature approaches dissipation in infinite-dimensional
systems from several complementary perspectives. Passive system theory
emphasizes input-output behavior and storage functions
\cite{BthK,Staf05}, while port-Hamiltonian formulations describe
partial differential equations directly in terms of energy and
interconnection structures \cite{JacZwa12,PhilippReisSchaller2025}.
We take a complementary viewpoint: starting from a given semigroup
generator, we ask whether its dynamics admit a dissipative Hamiltonian
representation and which properties the energy operator
must satisfy. Related work in the finite-dimensional
setting includes \cite{AchAC21,BeaMXZ26}; see
\cite{ArnEMN26} for preliminary results in infinite
dimensions.

More precisely, we consider the linear evolution equation
\begin{equation}\label{eq:cauchy}
 \dot x(t)=Ax(t),
\end{equation}
where $A:\dom(A)\subset\X\to\X$ generates a strongly continuous
semigroup $T=(T(t))_{t\geq0}$ on a Hilbert space $\X$.
Throughout the paper, Hilbert spaces may be real or complex, and all spaces
appearing in the same statement are understood to be over the same
scalar field.
We denote the norm and inner product on $\X$ by
$\norm{\cdot}_{\X}$ and $\inner{\cdot}{\cdot}_{\X}$,
respectively, and omit the subscripts whenever
the underlying space is clear.
We investigate when its dynamics can be represented in either of the two
forms
\begin{align}
 \dot x(t)&=D_{\mathrm e}Hx(t),\label{expevoleq}\\
 M\dot z(t)&=D_{\mathrm i}z(t),\label{impevoleq}
\end{align}
with dissipative operators $D_{\mathrm e}$ and $D_{\mathrm i}$ and
positive selfadjoint energy operators $H$ or $M$, respectively.
The dissipative operators $D_{\mathrm e}$ and $D_{\mathrm i}$
need not coincide. Throughout the paper, we assume that
$\ker H=\ker M=\{0\}$, but do not require either energy operator
to be coercive. 

The implicit formulation \eqref{impevoleq} is reminiscent
of the mass-operator formulations arising in weak formulations
of partial differential equations and their finite element
discretizations. We therefore also refer to it as a
\emph{mass representation}.

The case of energy operators with nontrivial kernels
is also of interest, but lies outside the scope
of the present paper; see \cite{ZwaM24}. For sufficiently regular solutions, the quadratic energies
of the explicit and implicit formulations satisfy, respectively,
\begin{align*}
 \frac{\mathrm d}{\mathrm dt}
 \frac12\norm{H^{1/2}x(t)}_{\X}^2
 &=\Real\inner{D_{\mathrm e}Hx(t)}{Hx(t)}_{\X}\leq0,\\
 \frac{\mathrm d}{\mathrm dt}
 \frac12\norm{M^{1/2}z(t)}_{\X}^2
 &=\Real\inner{D_{\mathrm i}z(t)}{z(t)}_{\X}\leq0.
\end{align*}
In infinite dimensions, however, these formal energy balances
do not automatically lead to well-defined operator representations.
For the explicit equation \eqref{expevoleq}, the product $D_{\mathrm e}H$ is defined
only for states $x\in\dom(H)$ with $Hx\in\dom(D_{\mathrm e})$.
Thus, recovering the original evolution requires not merely a
formal factorization, but equality of the operators $A$ and
$D_{\mathrm e}H$, including their domains.

The implicit formulation \eqref{impevoleq} involves a different difficulty.
If $M$ is unbounded, then the term $M\dot z$ is meaningful only when
$\dot z\in\dom(M)$, which need not hold for every classical
solution of the original evolution equation. Even for bounded
and injective $M$, the dissipative operator $D_{\mathrm i}=MA$
need not be closed. If $M$ is not surjective,
then $\lambda M-D_{\mathrm i}$ is not surjective
for any $\lambda$.

A further distinction arises when the energy is defined by a
closed quadratic form. Its dissipation inequality then may yield
contractive dynamics on a suitable energy space without
providing an operator factorization on the original state
space $\X$. We therefore distinguish exact operator
representations, implicit equations with specified solution
classes, and representations at the level of quadratic forms.
In particular, we do not assume a decomposition of the
dissipative operators into skew-adjoint and negative
semidefinite selfadjoint parts. For unbounded operators,
such a decomposition requires additional domain conditions.

A useful perspective comes from the renorming theory of semigroups.
Every uniformly bounded semigroup on a Banach space is
contractive with respect to the equivalent norm
\[
 \norm{x}_{*}:=\sup_{t\geq0}\norm{T(t)x}.
\]
More generally, every strongly continuous semigroup becomes quasi-contractive after
an equivalent Banach-space renorming; see \cite{EngeNage00}.
Even when the original space is a Hilbert space, such a norm need not be induced
by an inner product. In fact, there exist exponentially stable
Hilbert-space semigroups that are not similar, through a boundedly
invertible operator, to a contraction semigroup; see
\cite{LeMerdy2000,OlivaMazaTomilov2025}.
By contrast, a contractive \emph{equivalent Hilbertian} norm is induced
by a bounded, coercive positive operator. This observation motivates
allowing positive energy operators that are non-coercive or unbounded,
while retaining a Hilbert-space framework.

A natural way of constructing such energy operators is through \emph{intertwining
operators}. Suppose that a closed, densely defined, injective operator
$E:\dom(E)\subset\X\to\Z$ intertwines the strongly continuous semigroup $T$ with a contraction
semigroup $S$ on a Hilbert space $\Z$, in the sense that $\dom(E)$ is
$T$-invariant and
\begin{equation}\label{eq:intertwining-intro}
 ET(t)x=S(t)Ex\qquad \text{for all }x\in\dom(E),\quad t\geq0.
\end{equation}
The positive quadratic form $\formh[x]=\norm{Ex}_{\Z}^2$ is then
closed and nonincreasing along trajectories in its domain, and its
associated positive selfadjoint operator is $H=E^*E$.
Neither $E$ nor its inverse on $\ran E$ is required to be bounded.
This form-theoretic construction is closely related to existing work
on unbounded storage operators and pseudo-similarity, in particular
\cite{BthK}. It does \emph{not}, without further domain assumptions,
provide an operator identity of the form $A=DH$ or an equivalent
implicit equation in the original state space. Understanding this
additional step is central to the representation problem that we discuss in this paper.

A particularly useful realization of the intertwiner is furnished by
output trajectories. For a suitable observation operator $C$ with
values in a Hilbert space $\Y$, consider
\[
 (\Oop_Cx)(t)=CT(t)x\qquad\text{for all } t\geq0,
\]
where the domain consists of the initial states with square-integrable
output trajectories. This mapping intertwines $T$ with the contractive
left-shift semigroup on $L^2(0,\infty;\Y)$. The shift approach has a substantial history
\cite{Levan1990,Weiss1989Admissible}, as do admissibility and
observability criteria
\cite{GaoHou1999,GrabowskiCallier1996,HaakOuhabaz2012,JacobPartington2004Survey,
PartingtonPott2005,TucsnakWeiss2009}.
Closed, potentially unbounded output-trajectory operators are also
treated in \cite{BthK}. In particular, the existence of an
infinite-time admissible and exactly observable output operator is
characterized by strong stability together with similarity to a
contraction semigroup \cite[Theorem~8.11]{OlivaMazaTomilov2025}.
Here we use observation maps primarily as a \emph{construction tool}:
we allow one-sided estimates, or neither upper nor lower estimates,
and seek the resulting dissipative Hamiltonian representations rather
than new admissibility criteria.

Our central result is a fixed-energy representation
criterion: under the domain condition
$\dom(A)\subset\dom(H)$, the Lyapunov dissipation
inequality is equivalent to an exact factorization
$A=DH$ with a maximally dissipative operator $D$.
In particular, maximal dissipativity of the factor
can always be achieved by extension, without any
additional range condition. For bounded $H$, the
criterion is also equivalent to contractivity of the
associated quadratic energy and yields an implicit
representation. In particular, every
exponentially stable semigroup has both representations
with bounded energy. We then use pseudo-similarity and
output-trajectory maps to construct energies that need not be
bounded or coercive. For unbounded energy operators, a
contractive quadratic form always determines a contraction
semigroup on a canonical energy space, but a strong factorization
on the original space demands additional domain compatibility.
A dual construction and a variational interpretation complement
the operator result. We stress that the aim of the paper is a unified representation framework,
not a new general observability theorem.

The paper is organized as follows. \Cref{sec:prelim}
introduces the representation concepts and pseudo-similarity.
\Cref{sec:bounded} gives the fixed-energy criterion
and bounded-energy consequences. \Cref{sec:pseudo} treats closed intertwiners
and the output-trajectory construction. \Cref{sec:unbounded}
discusses energy-space dynamics, duality, and variational
realizations. \Cref{sec:examples} provides examples and
counterexamples, and \Cref{sec:conclusions} a conclusion.

\section{Preliminaries and representations}\label{sec:prelim}

Throughout the paper, Hilbert-space inner products are linear in the first
argument. All operators are linear, and equalities of operators are
understood to include equality of their domains. For a Hilbert space $\X$,
$\LL(\X)$ denotes the bounded operators on $\X$.
If $A$ generates
a strongly continuous semigroup $T$, then $A$ is densely defined and closed, and
$T(t)\dom(A)\subset\dom(A)$ with $AT(t)x=T(t)Ax$ for $x\in\dom(A)$;
see \cite[Lemma~II.1.3 and Theorem~II.1.4]{EngeNage00}.
A subspace $\mathcal D\subset\X\times\X$
is called \emph{dissipative} if
\[
 \Real\inner{u}{v}_{\X}\leq0\qquad\text{for all }(u,v)\in\mathcal D.
\]
An operator is dissipative if its graph is dissipative, and is
\emph{maximally dissipative} if its graph has no proper dissipative
linear extension. In particular, if a dissipative operator $D$ satisfies
$\ran(\lambda I-D)=\X$ for some $\lambda>0$, then
it is maximally dissipative by the injectivity property
for dissipative extensions
\cite[Proposition~II.3.14(i)]{EngeNage00}.
By the Lumer--Phillips theorem
\cite[Corollary~II.3.20]{EngeNage00},
it also generates a strongly continuous contraction semigroup;
see also \cite{LumerPhillips1961}. More general linear relations are useful for
singular Hamiltonian descriptions \cite{GerHaRe21}.

An operator $H$ will be called \emph{positive definite} if it is
selfadjoint, $H\geq0$, and $\ker H=\{0\}$; note that this does not mean that
$H\geq cI$ for some $c>0$. The latter property will be called
\emph{coercivity}. If $H$ is unbounded, its positive quadratic form
has domain $\dom(H^{1/2})$.

\begin{definition}\label{def:representation}
Let $A$ be the generator of a strongly continuous semigroup $T$ on a Hilbert space $\X$.
\begin{enumerate}[labelindent=1em,leftmargin=*,label={\rm(\roman*)}]
\item\label{itm:repr-explicit} An \emph{explicit dissipative Hamiltonian representation}
of $A$ is a factorization $A=DH$, where $H$ is a positive definite
selfadjoint operator and $D$ is dissipative.
\item\label{itm:repr-implicit} An \emph{implicit dissipative Hamiltonian representation} of
\eqref{eq:cauchy} consists of a positive definite selfadjoint operator
$M$ and a dissipative operator $D$ such that the equation
$M\dot z=Dz$ is meaningful and equivalent to \eqref{eq:cauchy}
for a specified class of solutions. If $M$ is bounded, a sufficient
condition for equivalence of classical solutions is $D=MA$ with
$\dom(D)=\dom(A)$.
\end{enumerate}
\end{definition}

Note that the specification of the solution class in
\ref{itm:repr-implicit} is essential. In particular,
writing $D=MA$ when $M$ is unbounded does not say that $MA$ is defined
on $\dom(A)$. The same care is required for the product $DH$.

\begin{definition}\label{def:pseudo}
Let $T_j$ be strongly continuous semigroups on Hilbert spaces $\X_j$ for $j=1,2$.
A \emph{pseudo-similarity from $T_1$ to $T_2$} is a closed, densely
defined, injective operator
$E:\dom(E)\subset\X_1\to\X_2$ such that $\dom(E)$ is
$T_1$-invariant and
\[
 E T_1(t)x=T_2(t)Ex\qquad \text{for all }x\in\dom(E),\quad t\geq0.
\]
We may, and often will, replace $\X_2$ by
$\overline{\ran E}$, which is invariant under $T_2$.
\end{definition}

The word pseudo-similarity refers here to a directed intertwining
relation; no boundedness of $E$ or $E^{-1}:\ran E\to\X_1$ is
assumed. Various concepts of quasi- and pseudo-similarity are studied
in operator theory; see, e.g., \cite[Definition~1.8]{BthK} and
\cite{OlivaMazaTomilov2025}.
\section{A fixed-energy representation criterion}\label{sec:bounded}
We begin by studying the operator-theoretic core of the
representation problem. Although the Hilbert-space Lyapunov inequality is classical, see e.g. 
\cite{EngeNage00,GaoHou1999,TucsnakWeiss2009}, the purpose of the following construction is to identify the corresponding dissipative Hamiltonian factors,
including their domains, within a single statement.

\begin{theorem}[Fixed-energy representation criterion]
\label{thm:bounded-factorizations}
Let $A$ be the generator of a strongly continuous semigroup $T$
on a Hilbert space $\X$, and let $H$ be a positive definite
selfadjoint operator on $\X$. Suppose that
\begin{equation}\label{eq:domain-condition}
 \dom(A)\subset\dom(H).
\end{equation}
Then the following assertions are equivalent:
\begin{enumerate}[labelindent=1em,leftmargin=*,
                  label={\rm(\alph*)}]
\item\label{itm:crit-lyap}
$\Real\inner{Ax}{Hx}\leq0$ for every $x\in\dom(A)$;
\item\label{itm:crit-factor}
there exists a dissipative operator $D$ on $\X$
such that $A=DH$, including equality of domains;
\item\label{itm:crit-maximal}
there exists a maximally dissipative operator $D$
on $\X$ such that $A=DH$, including equality of domains.
\end{enumerate}
If these conditions hold, the canonical dissipative factor
\begin{equation}\label{eq:d-explicit}
 \dom(D_{\mathrm e})=H\dom(A),\qquad
 D_{\mathrm e}(Hx)=Ax\quad\text{for all }x\in\dom(A)
\end{equation}
is densely defined. Every dissipative extension of
$D_{\mathrm e}$ preserves the factorization $A=DH$.
In particular, $D_{\mathrm e}$ admits a maximally
dissipative extension.

If $H$ is bounded, then the above assertions are also
equivalent to
\begin{equation}\label{eq:semigroup-H}
 T(t)^*HT(t)\leq H\qquad\text{for all }t\geq0.
\end{equation}
In this case the operator $D_{\mathrm i}=HA$ with
$\dom(D_{\mathrm i})=\dom(A)$ is dissipative, and the two equations
\begin{equation}\label{eq:factor-equivalence}
 \dot x=D_{\mathrm e}Hx,\qquad
 H\dot z=D_{\mathrm i}z
\end{equation}
both describe the dynamics generated by $A$ for
classical solutions.
Moreover, $E=H^{1/2}$ intertwines $T$ with a contraction
semigroup on $\overline{\ran E}=\X$.
\end{theorem}
\begin{proof}
Assume \ref{itm:crit-lyap}. Injectivity of $H$ makes
\eqref{eq:d-explicit} well-defined, and
\[
 \Real\inner{D_{\mathrm e}Hx}{Hx}
 =\Real\inner{Ax}{Hx}\leq0
 \qquad\text{for all }x\in\dom(A).
\]
Thus, $D_{\mathrm e}$ is dissipative.
Moreover, $D_{\mathrm e}H=A$, including equality
of domains, by the injectivity of $H$.

We first show that $D_{\mathrm e}$ is densely defined.
Choose $\mu>0$ such that $\mu I-A$ is bijective.
Let $z\in(H\dom(A))^\perp$ and set
$x=(\mu I-A)^{-1}z$. Then
\[
 0=\Real\inner{z}{Hx}
  =\mu\inner{x}{Hx}-\Real\inner{Ax}{Hx}
  \geq\mu\inner{x}{Hx}.
\]
Positive definiteness of $H$ implies that $x=0$,
and hence $z=0$. Therefore
$\overline{H\dom(A)}=\X$.

Since $D_{\mathrm e}$ is densely defined and
dissipative, its closure $D_0=\overline{D_{\mathrm e}}$
is dissipative. The range
$\mathcal R=\ran(I-D_0)$ is closed, since $D_0$ is closed
and
\[
 \norm{(I-D_0)u}\geq\norm{u}
\qquad\text{for all }u\in\dom(D_0).\]
Define an extension $D$ by
\[
 \dom(D)=\dom(D_0)\dotplus \mathcal R^\perp,\qquad
 D(u+v)=D_0u-v,
\]
where $u\in\dom(D_0)$ and $v\in \mathcal R^\perp$.
The sum is direct by the dissipativity of $D_0$.
Since $(I-D_0)u\perp v$, we obtain
\[
 \Real\inner{D(u+v)}{u+v}
 =\Real\inner{D_0u}{u}-\norm{v}^2\leq0.
\]
Furthermore,
\[
 \ran(I-D)=\mathcal R\oplus \mathcal R^\perp=\X.
\]
Thus $D$ is maximally dissipative and generates
a strongly continuous contraction semigroup.

It remains to show that every dissipative extension
$\widetilde D$ of $D_{\mathrm e}$ satisfies
$A=\widetilde D H$.
Since $\widetilde D$ extends $D_{\mathrm e}$,
we have $A\subseteq\widetilde D H$.
Let $x\in\dom(\widetilde D H)$.
By the surjectivity of $\mu I-A$, there exists
$y\in\dom(A)$ such that
\[
 (\mu I-A)y=(\mu I-\widetilde D H)x.
\]
With $w=x-y$ we obtain
$w\in\dom(\widetilde D H)$ and
$\widetilde D Hw=\mu w$. The dissipativity  then yields that
\[
 0\geq\Real\inner{\widetilde D Hw}{Hw}
   =\mu\inner{w}{Hw},
\]
and the positive definiteness of $H$ implies that $w=0$.
Consequently, $x=y\in\dom(A)$, which proves that
$A=\widetilde D H$.

We have shown that \ref{itm:crit-lyap} implies
\ref{itm:crit-maximal}.
The implication \ref{itm:crit-maximal}
$\Rightarrow$ \ref{itm:crit-factor} is immediate,
and \ref{itm:crit-factor} implies
\ref{itm:crit-lyap} by dissipativity.

Finally, suppose that $H$ is bounded.
Differentiating $\inner{HT(t)x}{T(t)x}$
for $x\in\dom(A)$, and using density,
proves the equivalence of \ref{itm:crit-lyap}
and \eqref{eq:semigroup-H}.
The identity $D_{\mathrm i}=HA$ is dissipative,
and injectivity of $H$ gives the stated classical
equivalence. On $\ran H^{1/2}$ define
$S(t)H^{1/2}x=H^{1/2}T(t)x$.
The inequality \eqref{eq:semigroup-H} then makes $S(t)$
contractive. It extends by density to a
strongly continuous contraction semigroup;
strong continuity follows from the boundedness
of $H^{1/2}$ and the strong continuity of $T$.
\end{proof}
\begin{remark}[Maximality of the canonical factor]
\label{rem:canonical-maximality}
For every $\lambda>0$ one has
\[
 \ran(\lambda I-D_{\mathrm e})
 =\ran(\lambda H-A).
\]
Therefore, the canonical factor $D_{\mathrm e}$
is maximally dissipative if and only if
$\ran(\lambda H-A)=\X$ for some $\lambda>0$.

If $A$ is bijective, then $HA^{-1}$ is closed
and everywhere defined, hence bounded. Therefore
$\lambda H-A=(\lambda HA^{-1}-I)A$ is bijective
for all sufficiently small $\lambda>0$.
In this case $D_{\mathrm e}$ itself is
maximally dissipative.

The range condition is not automatic, even when both
$A$ and $H$ are bounded. The second construction in
\Cref{subsec:output-no-factorization} shows that the
canonical factor may have a proper dense, nonclosed
domain, while its maximally dissipative extension
still yields the same exact factorization.
\end{remark}
Our next result discusses the implicit formulation.
\begin{proposition}[The mass-range resolvent]
\label{prop:mass-resolvent}
Let $M\in\LL(\X)$ be selfadjoint, positive definite, and suppose that
$\Real\inner{MAx}{x}\le0$ for $x\in\dom(A)$. Set $D_{\mathrm i}=MA$ with
$\dom(D_{\mathrm i})=\dom(A)$.
Then $D_{\mathrm i}$ is dissipative, and for every scalar
$\lambda$ such that $\lambda I-A$ is bijective, the pencil
\[
 L(\lambda)=\lambda M-D_{\mathrm i}:
 \dom(A)\longrightarrow\X
\]
is injective, with $\ran L(\lambda)=\ran M$. Thus $L(\lambda)$ is bijective as a map from
$\dom(A)$ onto $\ran M$, and
\begin{equation}\label{eq:relative-resolvent}
 L(\lambda)^{-1}M=(\lambda I-A)^{-1}
\quad\text{on }\X.
\end{equation}
In particular, the implicit equation \eqref{impevoleq} describes the dynamics
generated by $A$ for classical solutions, although, when
$M$ is not surjective, no
$L(\lambda)$ is surjective onto $\X$.
\end{proposition}
\begin{proof}
The assumed inequality $\Real\inner{MAx}{x}\leq0$
for $x\in\dom(A)$ shows that $D_{\mathrm i}$ is dissipative.
Fix $\lambda$ such that $\lambda I-A$ is bijective.
Then $L(\lambda)x=M(\lambda I-A)x$ for every $x\in\dom(A)$.
Since $M$ is injective and $\lambda I-A$ is bijective, this
identity proves injectivity, the range formula, and
\eqref{eq:relative-resolvent}. Classical equivalence follows
from the injectivity of $M$.
\end{proof}

\begin{remark}\label{rem:bounded-notm}
The implicit factor $D_{\mathrm i}=HA$ need not be closed,
let alone maximally dissipative, even though the two equations
in \eqref{eq:factor-equivalence} are classically equivalent.
\Cref{prop:mass-resolvent} shows why invertibility of the
pencil $L(\lambda)$ on all of $\X$ is not the right criterion for the dynamics
in this setting.
For unbounded $H$, condition \eqref{eq:domain-condition} is
an actual restriction, not a technical convenience.
The first construction in
\Cref{subsec:output-no-factorization} exhibits a
contractive closed energy form whose associated
energy operator does not satisfy this domain condition.
If the condition holds, form contractivity implies
\ref{itm:crit-lyap} in
\Cref{thm:bounded-factorizations}; see \Cref{sec:unbounded}.
\end{remark}
As a corollary we obtain that every exponentially stable semigroup admits both dissipative Hamiltonian representations.
\begin{corollary}\label{cor:expstable}
Every exponentially stable strongly continuous semigroup
on a Hilbert space admits both an explicit and an implicit
dissipative Hamiltonian representation with the same bounded,
positive definite energy operator $H$. For this choice of $H$,
the canonical explicit factor $D_{\mathrm e}$ is maximally
dissipative.
\end{corollary}
\begin{proof}
For the infinite-time Gramian, understood as a weak operator integral,
\begin{equation}\label{eq:gramian-id}
 H=\int_0^\infty T(t)^*T(t)\,dt,
\end{equation}
exponential stability implies boundedness, since
$\norm{T(t)}\le K e^{-\omega t}$ gives
$\norm{H}\le K^2/(2\omega)$, and
$\inner{Hx}{x}=\int_0^\infty\norm{T(t)x}^2dt>0$ for $x\ne0$.
Differentiating the energy for $x\in\dom(A)$ gives
$2\Real\inner{HAx}{x}=-\norm{x}^2$.
The integral representation in
\cite[Theorem~II.1.10]{EngeNage00} yields
$A^{-1}=-\int_0^\infty T(t)\,dt$.
Hence \Cref{thm:bounded-factorizations}
and \Cref{rem:canonical-maximality} apply.
\end{proof}

The construction of the Gramian \eqref{eq:gramian-id}
and its connection to Lyapunov equations are described
in \cite[Theorem~5.1.1]{TucsnakWeiss2009}. Conversely, if the identity
observation has square-integrable trajectories for all
$x\in\X$, the Datko--Pazy theorem yields exponential
stability; see \cite{Datko1970} and
\cite[Theorem~V.1.8]{EngeNage00}.

\begin{remark}[Lyapunov equations and dissipative Hamiltonian representations]
\label{rem:lyapunov}
For bounded $H$, the dissipation condition in
\Cref{thm:bounded-factorizations} is precisely the
Lyapunov inequality $A^*H+HA\leq0$, understood in the
quadratic-form sense on $\dom(A)$. For the Gramian \eqref{eq:gramian-id} of an exponentially
stable semigroup, the  following stronger operator identity holds:
$H\dom(A)\subset\dom(A^*)$ and
\[
 (A^*H+HA)x=-x
 \qquad\text{for all }x\in\dom(A).
\]
Indeed, polarization of the Gramian identity yields
\[
 \inner{HAx}{y}+\inner{Hx}{Ay}=-\inner{x}{y}
 \qquad\text{for all }x,y\in\dom(A),
\]
which proves the assertion.

In finite dimensions, positive definiteness of $H$
implies coercivity, and hence the existence of an
equivalent contractive Hilbertian norm. In infinite
dimensions, however, even the strict Lyapunov equation
need not yield a coercive solution. Thus exponential
stability does not generally imply similarity to a
contraction semigroup, as illustrated by the examples
in \Cref{sec:examples}.

The significance of \Cref{thm:bounded-factorizations}
is that even a non-coercive Lyapunov operator can yield
an exact factorization $A=DH$ with a maximally
dissipative factor $D$, while the canonical factor
$D_{\mathrm e}$ need not itself be maximally dissipative;
see \Cref{rem:canonical-maximality}.
\end{remark}

\section{Pseudo-similarity and output trajectories}\label{sec:pseudo}

\subsection{Closed energy forms from intertwiners}

\begin{theorem}\label{thm:pseudo-energy}
Let $\X$ and $\Z$ be Hilbert spaces, let
$E:\dom(E)\subset\X\to\Z$ be a closed, densely defined,
injective operator, and let $S$ be a contraction semigroup
on $\Z$.
Assume that $E$ intertwines $T$ with $S$ as in
\eqref{eq:intertwining-intro}. Then
\begin{equation}\label{eq:form}
 \formh(x,y)=\inner{Ex}{Ey}_{\Z},
 \qquad \dom(\formh)=\dom(E),
\end{equation}
is a densely defined, closed, positive definite form, satisfying
\begin{equation}\label{eq:formcontraction}
 \formh[T(t)x]\leq\formh[x]
 \qquad\text{for all } x\in\dom(E),\quad t\geq0.
\end{equation}
The associated positive selfadjoint operator is $H=E^*E$, with
\[
 \dom(H^{1/2})=\dom(E),\qquad
 \norm{H^{1/2}x}=\norm{Ex}.
\]
In particular, the semigroup $T$ is contractive on the completion
of $\dom(E)$ in the $\formh^{1/2}$-norm, which can be identified
with $\overline{\ran E}$.
\end{theorem}
\begin{proof}
The form in \eqref{eq:form} is closed because $E$ is closed;
it is strictly positive because $E$ is injective. The form
representation theorem for positive selfadjoint operators yields that
$H=E^*E$ and the asserted square-root identity; see
\cite[Theorems~VI.2.1 and~VI.2.23]{Kato80}.
By intertwining and the contractivity of $S$, we have
\[
 \norm{ET(t)x}=\norm{S(t)Ex}\leq\norm{Ex}.
\]
Finally the map $x\mapsto Ex$ identifies the completion in the
energy norm with $\overline{\ran E}$, and the induced semigroup
is precisely $S$ restricted to that closed invariant subspace.
\end{proof}

The converse to the construction in \Cref{thm:pseudo-energy} is also useful: a contractive
closed form produces an intertwining operator without an output
being specified beforehand.

\begin{proposition}\label{prop:closedform-converse}
Let $\formh$ be a densely defined, closed, positive definite,
nonnegative quadratic form on a Hilbert space $\X$, whose domain is $T$-invariant
and satisfies \eqref{eq:formcontraction}. Let $H$ be its associated
positive selfadjoint operator. Then $E=H^{1/2}$ is a closed,
densely defined, injective pseudo-similarity from $T$ to a
contraction semigroup on $\overline{\ran H^{1/2}}=\X$.
\end{proposition}
\begin{proof}
The map $S(t)$ defined initially by
$S(t)H^{1/2}x=H^{1/2}T(t)x$ is well-defined and contractive on
the dense subspace $\ran H^{1/2}$; it extends uniquely to a
contraction semigroup on $\X$. It remains to check strong
continuity. Fix $x\in\dom(H^{1/2})$ and let $t_n\downarrow0$.
The sequence $H^{1/2}T(t_n)x$ is bounded by the form inequality \eqref{eq:formcontraction}.
Every weakly convergent subsequence has weak limit $H^{1/2}x$:
indeed $T(t_n)x\to x$ in $\X$ and the graph of the closed
operator $H^{1/2}$ is weakly closed. Moreover,
\[
 \norm{H^{1/2}T(t_n)x}\leq\norm{H^{1/2}x}.
\]
Weak lower semicontinuity forces convergence of these norms,
and hence strong convergence to $H^{1/2}x$. Density and
contractivity give strong continuity on $\X$.
\end{proof}

The analysis in this subsection shows why closed energy forms provide a natural
language for pseudo-similar contraction dynamics. However, they do not,
by themselves, produce an everywhere valid operator product $DH$. We discuss this issue  in \Cref{sec:unbounded}.

\subsection{The output-trajectory transformation}

An observation operator $C:\dom(A)\to\Y$ that maps to a Hilbert space $\mathcal Y$ is assumed to be
bounded with respect to the graph norm
$\norm{x}_{\dom(A)}^2=\norm{x}_{\X}^2+\norm{Ax}_{\X}^2$.
It is called \emph{finite-time admissible} if there exist
$\tau>0$ and $\beta_\tau>0$ such that
\[
 \int_0^\tau\norm{CT(t)x}_{\Y}^2\,dt
 \leq\beta_\tau\norm{x}_{\X}^2
 \qquad\text{for all }x\in\dom(A).
\]
It is \emph{infinite-time admissible} if the same estimate
holds with $\tau=\infty$ and a finite constant $\beta_\tau$. An
infinite-time admissible observation operator $C$ is
called \emph{infinite-time exactly observable} if there
exists $\alpha>0$ such that
\[
 \int_0^\infty\norm{CT(t)x}_{\Y}^2\,dt
 \geq\alpha\norm{x}_{\X}^2
 \qquad\text{for all }x\in\dom(A).
\]
See \cite{PartingtonPott2005,TucsnakWeiss2009} for
these notions and related results.

We first consider $C\in\LL(\X,\Y)$, for which finite-time
admissibility is automatic. Define its possibly unbounded
infinite-time output operator by 
 \begin{align*}  \dom(\Oop_C)&=\left\{x\in\X:
 \int_0^\infty\norm{CT(t)x}_{\Y}^2\,dt<\infty\right\},\\
 (\Oop_Cx)(t)&=CT(t)x\quad\text{for all }t\geq0.
 \end{align*}
The left-shift semigroup on $L^2(0,\infty;\Y)$ is
\[
 (S(t)f)(\tau)=f(t+\tau)\qquad\text{for all }t\geq0
\text{ and a.e.\ }\tau\geq0.
\]
It is a strongly continuous contraction semigroup and it is strongly
stable.  The next proposition records elementary properties which
also motivate more general output maps for graph-bounded $C$;
see \cite{Levan1990},
\cite[Proposition~3.1]{BthK}, and
\cite[Secs.~4.3 and~4.6]{TucsnakWeiss2009}.
\begin{proposition}\label{prop:output-closed}
For $C\in\LL(\X,\Y)$, $\Oop_C$ is closed, its domain is
$T$-invariant, and
\[
 \Oop_CT(t)x=S(t)\Oop_Cx,
 \qquad\text{for all } x\in\dom(\Oop_C),\; t\geq0.
\]
If, in addition, $\dom(\Oop_C)$ is dense and
$CT(\cdot)x=0$ implies $x=0$ for $x\in\dom(\Oop_C)$,
then $\Oop_C$ is injective (the qualitative or approximate
observability condition) and defines a pseudo-similarity
to the shift semigroup.
Moreover, every $x\in\dom(\Oop_C)$ satisfies
\[
 \norm{\Oop_C T(t)x}_{L^2}^2-\norm{\Oop_Cx}_{L^2}^2
 =-\int_0^t\norm{CT(\tau)x}_{\Y}^2\,d\tau\qquad\text{for all }t\geq0.
\]
\end{proposition}
\begin{proof}
Suppose that $x_n\in\dom(\Oop_C)$, $x_n\to x$ in $\X$, and
$\Oop_Cx_n\to f$ in $L^2(0,\infty;\Y)$. On every finite
interval $[0,\tau]$, one has the convergence $CT(\cdot)x_n\to CT(\cdot)x$ uniformly.
Hence, $f=CT(\cdot)x$ almost everywhere and
$x\in\dom(\Oop_C)$, proving closedness. The intertwining
identity follows from the semigroup property, and invariance follows
since the tail of an $L^2$-function is in $L^2$. The last
identity is obtained by subtracting the integrals over
$[t,\infty)$ and $[0,\infty)$.
\end{proof}

The upper admissibility estimate and the lower observability
estimate have a particularly simple interpretation in terms
of $\Oop_C$. Assume that $\Oop_C$ is densely defined and
injective. Then
\[
 \begin{array}{rcl}
 \Oop_C\text{ bounded}&\Longleftrightarrow&
 \displaystyle\int_0^\infty\norm{CT(t)x}^2\,dt
      \leq \beta\norm{x}^2,\\[2mm]
 \Oop_C^{-1}\text{ bounded on }\ran\Oop_C&\Longleftrightarrow&
 \displaystyle\int_0^\infty\norm{CT(t)x}^2\,dt
      \geq\alpha\norm{x}^2,
 \end{array}
\]
where the respective estimates hold for some
$\alpha,\beta>0$ and all $x\in\dom(\Oop_C)$.
The upper estimate is equivalent to infinite-time
admissibility of $C$: by closedness and density it implies
$\dom(\Oop_C)=\X$. Under this assumption, the lower
estimate is precisely the infinite-time exact-observability
inequality for $(A,C)$.

Without admissibility, the lower estimate may still hold
on the possibly proper domain $\dom(\Oop_C)$. This is a lower observability estimate on the natural
domain of the output transformation, but does not imply
infinite-time admissibility. In this case
$\ran\Oop_C$ is closed. Thus, bounded invertibility on the range of the output
transformation does not require admissibility of the
underlying observation operator.

When both estimates hold, the associated semigroup is similar to
a contraction semigroup. A precise existence characterization,
including the required strong stability, is given in
\cite[Theorem~8.11]{OlivaMazaTomilov2025}; compare the earlier
results in \cite{BthK,GrabowskiCallier1996}. We quote these results
rather than reproving their observation-theoretic content.

The connection can also be formulated intrinsically in terms
of Lyapunov operators. A bounded injective output map induces a
bounded positive definite $H=\Oop_C^*\Oop_C$ such that
$T(t)^*HT(t)\le H$ and
$\norm{H^{1/2}T(t)x}\to0$ as $t\to\infty$.
Conversely, given such an $H$, the nonnegative form
\[
 (x,y)\longmapsto
 -\inner{HAx}{y}-\inner{Hx}{Ay}\qquad\text{for all } x,y\in\dom(A),
\]
is bounded with respect to the graph norm on $\dom(A)$.
Factoring this nonnegative form yields a graph-norm bounded
observation operator $C:\dom(A)\to\Y$, for a suitable
Hilbert space $\Y$, such that
\[
 \int_0^\infty\norm{CT(t)x}^2\,dt=\inner{Hx}{x}
 \qquad\text{for all }x\in\dom(A).
\]
In particular, $C$ is infinite-time admissible.
If $H$ is additionally coercive, the pair $(A,C)$
is infinite-time exactly observable.
These standard Lyapunov and admissibility constructions
are discussed in \cite{GaoHou1999,OlivaMazaTomilov2025};
see also \cite[Sec.~5.1]{TucsnakWeiss2009}.

\subsection{Both one-sided output constructions for stable semigroups}

There is a particularly simple way to construct the two one-sided
cases without asserting any similarity to a contraction semigroup.

\begin{theorem}\label{thm:two-outputs}
Let $T$ be exponentially stable. Then:
\begin{enumerate}[labelindent=1em,leftmargin=*,label={\rm(\alph*)}]
\item\label{itm:outputs-identity} The identity observation $C_0=I$ produces a bounded,
injective infinite-time output map
$\Oop_0x=T(\cdot)x$.
\item\label{itm:outputs-generator} The graph-norm bounded observation $C_1=A$, initially
defined by
\[
 \Oop_1x=T(\cdot)Ax\qquad\text{for all } x\in\dom(A),
\]
is closable as an operator from $\X$ into $L^2(0,\infty;\X)$.
Its closure $\overline{\Oop_1}$ is densely defined, closed,
injective, intertwines $T$ with the left shift, and has a bounded
inverse on its range.
\end{enumerate}
Consequently, every exponentially stable semigroup has both a
bounded injective output transformation and a closed, bounded-below
output transformation, even if it is not similar to a contraction
semigroup.
\end{theorem}
\begin{proof}
Choose $K\ge1$ and $\omega>0$ such that
$\norm{T(t)}\le K e^{-\omega t}$ for all $t\ge0$.
Part~\ref{itm:outputs-identity} follows immediately from
\[
 \int_0^\infty\norm{T(t)x}^2\,dt
 \leq \frac{K^2}{2\omega}\norm{x}^2.
\]
For $x\in\dom(A)$, commutation with the semigroup and exponential
stability show $\Oop_1x\in L^2$ and
\[
 \norm{\Oop_1x}_{L^2}
 \leq\frac{K}{\sqrt{2\omega}}\norm{Ax}.
\]
Moreover, differentiating the squared norm, integrating over
$[0,\infty)$, and applying the Cauchy--Schwarz inequality yield
\begin{align*}
 \norm{x}^2
 &=-2\Real\int_0^\infty
         \inner{AT(t)x}{T(t)x}\,dt\\
 &\leq2\norm{\Oop_1x}_{L^2}\norm{\Oop_0x}_{L^2}
 \leq2\norm{\Oop_0}\norm{\Oop_1x}_{L^2}\norm{x}.
\end{align*}
It follows that
\begin{equation}\label{eq:lowerC=A}
 \norm{\Oop_1x}_{L^2}\geq
 \frac{1}{2\norm{\Oop_0}}\norm{x}
 \qquad\text{for all } x\in\dom(A).
\end{equation}
To see closability, let $B$ be the generator of the left-shift
semigroup, $Bf=f'$ on $H^1(0,\infty;\X)$.
For $x\in\dom(A)$ we have
$\Oop_0x\in\dom(B)$ and $B\Oop_0x=\Oop_1x$.
Since $\Oop_0$ is bounded and $B$ is closed,
$B\Oop_0$ is closed. As $\Oop_1$ is its restriction
to $\dom(A)$, it is closable.
The estimate \eqref{eq:lowerC=A} extends to its closure,
proving
injectivity and boundedness of the inverse. The intertwining
identity holds on $\dom(A)$ and passes to the closure using
boundedness of $T(t)$ and $S(t)$.
\end{proof}

Although $C_1=A:\dom(A)\to\X$ is graph-norm bounded,
this does not ensure finite-time or infinite-time
admissibility. For exponentially stable semigroups,
these two notions of admissibility are equivalent;
see \cite[Sec.~4.6]{TucsnakWeiss2009}.
Its closed output transformation
$\overline{\Oop_1}$ is bounded below. If it were also bounded,
the associated semigroup would be similar to a contraction
semigroup. Hence, for an exponentially stable semigroup not
similar to a contraction, $\overline{\Oop_1}$ must be unbounded,
while $\Oop_0^{-1}$ is unbounded on its range.

The existence of a closed, bounded-below pseudo-similarity in
\Cref{thm:two-outputs} cannot be extended to arbitrary uniformly
bounded strongly continuous semigroups. Indeed, there exist 
semigroups $T$ for which no bounded, injective, dense-range
intertwiner $F$ from a strongly continuous contraction semigroup
$S$ satisfies $T(t)F=FS(t)$ for all $t\geq0$.
The norm bound of $T$ can even be chosen arbitrarily close to one
\cite[Theorem~1.4 and Remark~6.3]{OlivaMazaTomilov2026Tensor}.
In particular, such a semigroup cannot admit a closed, densely
defined, bounded-below pseudo-similarity to a contraction
semigroup: its inverse on the closed range would provide
precisely such an intertwiner.

\section{Unbounded energies and operator representations}
\label{sec:unbounded}

Let $E:\dom(E)\subset\X\to\Z$ be a closed, densely defined,
injective intertwiner from $T$ to a contraction semigroup $S$.
Write
\[
 \V=\dom(E),\quad \norm{x}_{\V}^2=\norm{x}_{\X}^2+\norm{Ex}_{\Z}^2,
 \quad \mathcal K=\overline{\ran E},\quad H=E^*E.
\]
The transition to the energy space draws on closed-form theory
\cite[Chap.~VI, Sec.~2]{Kato80} and the theory of
pseudo-similarity and storage functions
\cite[Sec.~9.2]{Staf05};
see also \cite[Secs.~1 and~7]{BthK}.
We record it to specify which operator is maximally dissipative.

\begin{proposition}[Energy-space realization and variational equation]
\label{thm:weakimplicit}
The restricted semigroup $T_{\V}=T|_{\V}$ is strongly continuous on
$\V$, with generator $A_{\V}$, the part of $A$ in $\V$:
\[
 \dom(A_{\V})=\{x\in\dom(A)\cap\V:Ax\in\V\},\qquad A_{\V}x=Ax.
\]
If $B$ generates $S|_{\mathcal K}$, then $B$ is maximally dissipative
on $\mathcal K$,
\[
 BEx=EA_{\V}x\quad\text{for all }x\in\dom(A_{\V}),
\]
and $E\dom(A_{\V})$ is a core for $B$.
Let $\V'$ denote the anti-dual of $\V$, that is,
the space of continuous conjugate-linear functionals
on $\V$. Define
\[
 \begin{aligned}
  \mathcal M &: \V\to\V',&
  \langle\mathcal Mx,y\rangle_{\V',\V}&=\inner{Ex}{Ey}_{\Z},\\
  \mathcal D &: \dom(A_{\V})\subset\V\to\V',&
  \mathcal D x&=\mathcal M A_{\V}x.
 \end{aligned}
\]
Then $\mathcal M$ is bounded and injective, $\mathcal D$ is
dissipative in the duality sense, and
\begin{equation}\label{eq:weak-implicit}
 \mathcal M\dot z=\mathcal D z\quad\text{in }\V'
\end{equation}
is equivalent to $\dot z=A_{\V}z$ for $\V$-classical solutions.
\end{proposition}
\begin{proof}
Intertwining gives
$\norm{T(t)x}_{\V}^2\le\norm{T(t)}^2\norm{x}^2+\norm{Ex}^2$.
Since $T(t)x\to x$ and $ET(t)x=S(t)Ex\to Ex$, the
restricted semigroup is strongly continuous on $\V$. The generator formula follows from
\cite[Proposition~2.4.4]{TucsnakWeiss2009}. Differentiating
$ET_{\V}(t)x=S(t)Ex$ gives $BEx=EA_{\V}x$.
For sufficiently large $\lambda$,
\[
 E(\lambda-A_{\V})^{-1}=(\lambda-B)^{-1}E
 \quad\text{on }\V.
\]
As $\ran E$ is dense in $\mathcal K$, the resolvent maps it onto
a graph-norm dense subset of $\dom(B)$, giving the core property.
Finally, $\mathcal M$ is bounded and injective since $E$ is,
and
\[
 2\Real\langle\mathcal D x,x\rangle_{\V',\V}
  =\left.\frac{\mathrm d}{\mathrm dt}\right|_{t=0+}
       \norm{S(t)Ex}_{\Z}^2\le0
\qquad\text{for all }x\in\dom(A_{\V}).\]
Injectivity of $\mathcal M$ yields the equivalence.
\end{proof}

There is also a useful \emph{dual} consequence of the lower
observation estimate. It concerns the adjoint dynamics and must not
be confused with an operator identity for $A$ on $\X$.

\begin{proposition}[Duality for bounded-below intertwiners]
\label{prop:dual-one-sided}
In addition to the hypotheses in~\Cref{thm:weakimplicit}, assume that
$\norm{Ex}\ge c\norm{x}$ on $\dom(E)$ for some $c>0$.
Then $H\ge c^2 I$ and $P=H^{-1}\in\LL(\X)$ is positive definite,
with
\[
 T(t)P T(t)^*\le P\qquad\text{for all } t\ge0.
\]
Consequently, the adjoint generator $A^*$ admits
a bounded-energy factorization $A^*=D_{\mathrm a}P$
with maximally dissipative $D_{\mathrm a}$.
\end{proposition}
\begin{proof}
Closedness and the lower bound imply $\ran E=\mathcal K$.
Thus $F=E^{-1}:\mathcal K\to\X$ is bounded, with dense range,
and $T(t)F=F S(t)|_{\mathcal K}$. Moreover $P=FF^*$, by the polar
decomposition of $E$. Therefore,
\[
 T(t)P T(t)^*=F S(t)S(t)^*F^*\le FF^*=P.
\]
The last assertion follows by applying~\Cref{thm:bounded-factorizations} to $A^*$ and $P$.
\end{proof}

\paragraph{From energy forms to operators on $\X$}
Contractivity of the energy form alone does not imply an
operator factorization on the original state space. Indeed,
if $H=E^*E$ and $x\in\dom(A)\cap\dom(H)$, then
\[
 0\geq\norm{H^{1/2}T(t)x}^2-\norm{H^{1/2}x}^2
 \geq 2\Real\inner{T(t)x-x}{Hx}.
\]
Dividing by $t$ and taking the limit to zero gives
$\Real\inner{Ax}{Hx}\leq0$. If
$\dom(A)\subset\dom(H)$, \Cref{thm:bounded-factorizations}
therefore provides an exact representation $A=DH$
with maximally dissipative $D$.
The domain inclusion is also necessary for an exact product
representation with this particular $H$, as the theorem shows.
The first construction in
\Cref{subsec:output-no-factorization} shows that it may
fail even for a closed, injective output-trajectory map.

\begin{remark}[Strong implicit equations and duality]
\label{rem:implicit-duality}
Suppose that $M$ is a closed operator and an $\X$-valued equation
$M\dot z=D_{\mathrm i}z$ is required to hold at $t=0$ for
\emph{every} classical trajectory of $\dot z=Az$. This requires
$\ran A\subset\dom(M)$. If $A$ is bijective, then
$\dom(M)=\X$, so the closed graph theorem forces $M$ to be
bounded. Thus, an unbounded mass operator requires a smaller solution
class or an equation in a different space, such as
\eqref{eq:weak-implicit}. The latter is not an operator-pair
representation on $\X$ in the sense of \Cref{def:representation}.
When $H$ is bounded and coercive, the change of variables
$w=Hx$ transforms $\dot x=D_{\mathrm e}Hx$ into
$H^{-1}\dot w=D_{\mathrm e}w$.
Without coercivity the change of variables need not be a
boundedly invertible coordinate transformation of $\X$.
\end{remark}

\section{Examples and limitations}\label{sec:examples}
The following examples illustrate the scope and limitations
of the representation results.
\subsection{All four boundedness combinations}

Let $I\subset(0,\infty)$ be an interval, set $\X=L^2(I)$, and let
\[
 (Ax)(s)=-s x(s),\qquad
 \dom(A)=\{x\in L^2(I):sx(s)\in L^2(I)\}.
\]
The semigroup $T(t)x(s)=e^{-st}x(s)$ is selfadjoint,
contractive, and strongly stable. For the identity observation,
\[
 \norm{\Oop_0x}^2_{L^2(0,\infty;\X)}
 =\frac12\int_I\frac{|x(s)|^2}{s}\,ds.
\]
The corresponding energy operator is multiplication by $1/(2s)$.
Consequently each of the following cases is realized:
\begin{center}
\begin{tabular}{lll}
\hline
$I$ & $\Oop_0$ & $\Oop_0^{-1}$ on $\ran\Oop_0$\\
\hline
$(1,2)$ & bounded & bounded\\
$(1,\infty)$ & bounded & unbounded\\
$(0,1)$ & unbounded & bounded\\
$(0,\infty)$ & unbounded & unbounded\\
\hline
\end{tabular}
\end{center}

Thus, the first case satisfies both the upper and lower
observation estimates, the second only the upper estimate,
the third only the lower estimate on its natural domain,
and the fourth neither. In particular, a lower observation
estimate does not imply admissibility.

In every case, the output map is closed, densely defined, and
injective. Notice that these four patterns can occur although the
original semigroup is already contractive. The distinction concerns
the particular \emph{energy representation}, not just the
contractivity of the semigroup itself.

\subsection{Limitations of canonical factorizations}
\label{subsec:output-no-factorization}

Let $\X=\ell^2(\mathbb N)$ and
$A=\operatorname{diag}(-1/n)$. Then $A$ is bounded and
$T(t)=\operatorname{diag}(e^{-t/n})$ is contractive and
strongly stable, but not exponentially stable. For $C=I$,
\[
 \norm{\Oop_Cx}_{L^2}^2=\frac12\sum_{n\ge1}n|x_n|^2,
 \qquad H=\Oop_C^*\Oop_C=\operatorname{diag}(n/2).
\]
The output map is closed, densely defined, injective, and
bounded below, but $H$ is unbounded. More precisely,
$\dom(\Oop_C)=\{x\in\ell^2:\sum_{n\ge1}n|x_n|^2<\infty\}$,
whereas
$\dom(H)=\{x\in\ell^2:\sum_{n\ge1}n^2|x_n|^2<\infty\}$. Since
$\dom(A)=\X\not\subset\dom(H)$, it cannot yield any exact
factorization $A=DH$ on $\X$, even though the energy is
nonincreasing. The obstruction concerns this particular $H$;
trivially, $A=A I$ is already a dissipative representation.
The same semigroup also illustrates the distinction
between the canonical dissipative factor and its
maximally dissipative extensions. Indeed,
$H_{\mathrm b}=\operatorname{diag}(1/n)$ is bounded,
positive definite, and a contractive energy, but
$D_{\mathrm e}=-I$ on $\ran H_{\mathrm b}$ is not
closed and hence not maximally dissipative.
Nevertheless, $D=-I$ on $\X$ is a maximally
dissipative extension of $D_{\mathrm e}$, and
$A=DH_{\mathrm b}$ holds with equality of domains.

\subsection{The Dirichlet heat equation}

Let $\Omega\subset\mathbb R^d$ be a bounded domain
and let $A=\Delta_D$ be the Dirichlet Laplacian
on $L^2(\Omega)$. Then $A$ is selfadjoint and
there exists $\lambda_1>0$ such that
\[
 \inner{Ax}{x}\leq-\lambda_1\norm{x}^2
 \qquad\text{for all }x\in\dom(A).
\]
Consequently, the corresponding heat semigroup
is exponentially stable.
By the functional calculus for selfadjoint operators
(see \cite[Chap.~VI, Sec.~5]{Kato80}), the identity
observation yields
\[
 H_0=-\frac12 A^{-1},
 \qquad D_{\mathrm e}=-2A^2
 \quad\text{on }\dom(A^2).
\]
The bounded operator $H_0$ is positive definite but not coercive.
The implicit factor is $D_{\mathrm i}=-\tfrac12 I$ restricted
to $\dom(A)$. With observation $C_1=A$, we instead obtain
\[
 H_1=-\frac12 A,
 \qquad D_{\mathrm e}=-2I\quad\text{on }\X.
\]
The closure of the output map initially defined by
$\Oop_1x=T(\cdot)Ax$ on $\dom(A)$ has domain
$\dom((-A)^{1/2})$, whereas $\dom(H_1)=\dom(A)$.
The energy operator $H_1$ is unbounded and coercive, but $H_1A=-\tfrac12 A^2$ is defined only on $\dom(A^2)$, not on
all of $\dom(A)$. On the other hand, for $M=H_0$ and
$D_{\mathrm i}=H_0A=-\tfrac12 I|_{\dom(A)}$, the
operator pencil $\lambda H_0-D_{\mathrm i}$ has range
$\ran H_0=\dom(A)\ne\X$
whenever $\lambda I-A$ is bijective, while \eqref{eq:relative-resolvent} recovers the heat resolvent.
Thus, both the explicit and implicit representations can be
meaningful without a boundedly invertible energy operator.

\subsection{Exponential stability without Hilbertian similarity}

Exponentially stable analytic Hilbert-space semigroups
need not be similar to contraction semigroups
\cite{LeMerdy1998,LeMerdy2000,OlivaMazaTomilov2025}.
More strikingly, there exist strongly continuous semigroups
that are nilpotent and immediately compact, as well as
strongly continuous semigroups that are quasi-nilpotent,
immediately compact, and bounded holomorphic of angle $\pi/2$.
Neither class is similar even to a quasi-contraction semigroup.
The construction uses a tensor-product splitting principle:
similarity of a tensor product to a quasi-contraction
semigroup forces the same property for its nonzero factors
\cite[Theorems~1.1 and~1.3]{OlivaMazaTomilov2026Tensor}.

Both classes of examples are exponentially stable.
\Cref{cor:expstable} therefore furnishes the bounded,
positive definite Gramian $H_0$ from \eqref{eq:gramian-id}
and an exact factorization $A=D_{\mathrm e}H_0$
with maximally dissipative $D_{\mathrm e}$.
For a nilpotent semigroup, the Gramian integral even
terminates at a finite extinction time, and in both cases
the strict identity
$2\Real\inner{H_0Ax}{x}=-\norm{x}^2$
holds for all $x\in\dom(A)$.
Nevertheless, $H_0$ cannot be coercive:
any bounded, coercive energy operator satisfying
\eqref{eq:semigroup-H} would yield a boundedly invertible
similarity to a contraction semigroup, contradicting
the examples. Thus strict Lyapunov dissipation,
even together with very strong stability and compactness,
does not imply an equivalent contractive Hilbertian norm.

The one-sided output constructions in
\Cref{thm:two-outputs} reflect the same obstruction:
the identity output is bounded with unbounded inverse
on its range, whereas the closure of the generator
output is unbounded but bounded below.

\section{Concluding remarks}\label{sec:conclusions}

We have characterized exact dissipative Hamiltonian
factorizations of semigroup generators in terms of
a Lyapunov inequality and a domain inclusion.
Whenever such a factorization exists, its dissipative
factor can be chosen maximally dissipative without
additional assumptions. The canonical factor need
not itself be maximally dissipative, but every
dissipative extension preserves the exact
operator factorization.
For bounded energy operators, the same Lyapunov inequality
also yields an implicit representation.

Closed intertwiners and output-trajectory maps provide
energy operators and contractive dynamics on associated
energy spaces. Bounded-below intertwiners additionally
yield a dual Lyapunov inequality and a dissipative Hamiltonian
factorization of the adjoint generator.
These constructions do not, in general, guarantee an exact
factorization of the original generator with the constructed
energy operator. Domain compatibility remains essential.

\section*{Acknowledgments}
Both authors gratefully acknowledge funding by the
Deutsche For\-schungsgemeinschaft (DFG, German Research
Foundation) -- Project-ID 531152215 -- CRC 1701
``Port-Hamiltonian Systems''.

\bibliographystyle{siamplain}
\bibliography{references}
\end{document}